\documentstyle[]{article}
\newcommand{\kitem}{\begin{itemize}\vspace{-2ex}}
\newcommand{\kenditem}{\vspace{-1ex}\end{itemize}}

\newcommand{\bean}{\[\begin{array}{rcl}}
\newcommand{\eean}{\end{array}\]}

\newcommand{\R}{{I\!\!R}}
\newcommand{\Z}{{Z\!\!\!Z}}
\newcommand{\Zz}{\Z}
\newcommand{\Hh}{I\!\!H}
\newcommand{\Fh}{I\!\!F}

\newcommand{\Rr}{\R}

\newcommand{\Ll}{I\!\!L}
\newcommand{\Pp}{I\!\!P}
\newcommand{\CO}{{\cal O}}
\newcommand{\CM}{{\cal M}}

\newcommand{\supp}{\mbox{\rm supp}\,}

\newcommand{\orb}{\mbox{\rm orb}\,}

\newcommand{\surj}{\longrightarrow\hspace{-1.5em}\longrightarrow}

\newcommand{\ko}{\overline}
\newcommand{\ku}{\underline}
\newcommand{\ks}{\scriptstyle}

\newcommand{\keps}{\varepsilon}
\newcommand{\veee}{\scriptscriptstyle\vee}

\newcommand{\gExt}{\mbox{\rm Ext}\,}

\newcommand{\gEnd}{\mbox{\rm End}\,}
\newcommand{\lEnd}{End\,}

\newcommand{\UB}{{\cal U}}

\newcommand{\kQ}{Q}
\newcommand{\wt}{\theta}
\newcommand{\cwt}{\theta^c}

\newcommand{\kTh}{{\cal T}(\wt)}

\newcommand{\kA}{{C}}

\newcommand{\kSwt}{{\Sigma(\wt)}}
\newcommand{\kScwt}{{\Sigma(\cwt)}}
\newcommand{\wtM}{{{\cal M}^{\wt}(\kQ)}}
\newcommand{\cwtM}{{{\cal M}(\kQ)}}

\newcommand{\kD}{\Delta}

\newcommand{\kMwt}{M^\wt}

\newcommand{\Alg}{\cA}

\newcounter{Abschnitt}[section]
\newcommand{\neu}[1]{\protect\refstepcounter{Abschnitt}\protect
   \label{#1}\vspace{1ex}
   {\bf (\protect\arabic{section}.\protect\arabic{Abschnitt})}
                     $\qquad$}
\newcommand{\zitat}[2]{(\protect\ref{#1}.\protect\ref{#1-#2})}

\newcommand{\DP}{.}

\newcommand{\mCoh}{\mbox{Coh}}

\newcommand{\mEnd}{\mbox{End}}
\newcommand{\mExt}{\mbox{Ext}}

\newcommand{\mforall}{\mbox{ for all }}

\newcommand{\mHom}{\mbox{Hom}}

\newcommand{\mmodstr}{\mbox{mod--}}

\newcommand{\mrad}{\mbox{rad}}

\newcommand{\ck}{{{k}}} 
 
\newcommand{\cA}{{\cal{A}}}

\newcommand{\cC}{{\cal{C}}}
\newcommand{\cD}{{\cal{D}}}
\newcommand{\cE}{{\cal{E}}}

\newcommand{\cL}{{\cal{L}}}
\newcommand{\cM}{{\cal{M}}}

\newcommand{\cR}{{\cal{R}}}

\newcommand{\cT}{{\cal{T}}}

\newcommand{\ra}{\rightarrow}
\newcommand{\lra}{\longrightarrow}

\begin{document}
\title{A Vanishing Result for the Universal Bundle on a Toric Quiver Variety}

\author{Klaus Altmann \qquad Lutz Hille}
\date{}
\maketitle

\begin{abstract} 
Let $Q$ be a finite quiver without oriented cycles. Denote
by $\UB \ra \cwtM$ the fine moduli space of stable thin sincere
representations of $Q$ with respect to the canonical stability notion. We 
prove $\mExt^l_{\cwtM}(\UB,\UB) = 0$ for all $l>0$ and compute the
endomorphism algebra  of the universal bundle $\UB$. Moreover, we obtain a
necessary and sufficient condition for when this algebra 
is isomorphic to the path algebra of the quiver
$Q$. If so, then the bounded derived categories of
finitely generated right $\ck Q$-modules and that  
of coherent sheaves on $\cwtM$ are related via the full and 
faithful functor $- \otimes^{\Ll}_{\ck Q}\UB$.
\end{abstract}


%
%
\section{Introduction}\label{Int}


\neu{Int-1}
Let $Q$ be a quiver (i.e.\ an oriented graph) without oriented cycles; 
denote by $Q_0$ the vertices and by $Q_1$ the arrows of $Q$. For a fixed
dimension vector $d$, that is a map $d: Q_0 \ra \Zz_{\geq 0}$, we
denote by $\Hh(d) := \{\wt: Q_0 \ra \Rr\mid \sum_{q \in Q_0} \wt_qd_q = 0
\}$ the vector space of the so-called weights with respect to $d$. 
We fix an algebraically closed field $k$. To each
$\wt\in\Hh(d)$ there exists the 
moduli space $\cM^{\wt}(Q,d)$ of $\wt$-semistable 
$k$-representations of $Q$ with dimension vector $d$ (cf.\ \cite{King}). 
This space is known to be projective
and, in case $\wt$ is in general position 
and $d$ is indivisible, also smooth. Moreover, if we
restrict ourselves to thin sincere representations, that is $d_q = 1$ for all $q \in
Q_0$, then $\wtM$ is also toric (cf.\ \cite{HilleTor}). 
In any case, each integral weight $\wt$
induces an ample line bundle $\cL(\wt)$ on $\cM^{\wt}(Q,d)$.\\ 
If $\wt$ is in general position and $d$ indivisible, then
$\cM^{\wt}(Q,d)$ is, in addition, a  
fine moduli space admitting a universal bundle $\UB$. The universal bundle
splits into a direct  
sum of vector bundles $\UB = \oplus_{q \in Q_0} \UB_q$, and the summands
$\UB_q$ have rank $d_q$ (cf.\ \cite{King}).
All known examples suggest that the universal bundles on those moduli
spaces have no self-extensions, i.e.\ 
$\gExt^{\, l}_{\cM^{\wt}(Q,d)}(\UB,\UB) = 0$ for all $l >0$. The issue of
this paper is to prove this formula in special cases.  The meaning of
this property and its relation to tilting theory will be discussed in
\zitat{Int}{4}. \\
In this paper we restrict ourselves to thin sincere representations; the
corresponding moduli spaces are called {\em toric quiver varieties}. Because
$d = (1,\ldots,1)$ is fixed, we will omit it in all notation introduced
above. The direct summands of the universal bundle are line
bundles, and they  
are characterized, up to a common twist, by the following property:
For any arrow $\alpha\in Q_1$ pointing from $p$ to $q$ ($p,q\in Q_0$) the
invertible sheaf $\,\UB_p^{-1} \otimes \UB_q$ corresponds to the
divisor of all representations assigning the zero map to $\alpha$.
Furthermore, there exists a distinguished weight $\cwt$ (see \zitat{Int}{2}
and \zitat{Mod}{7} for a definition and first properties). We denote the
corresponding moduli space by $\cwtM$. 
\par


\neu{Int-2}
Polarized projective toric varieties may be constructed from lattice
polytopes. If one wants to forget about the polarization, simply consider
the inner normal fan of the polytope. 
In \S \ref{Mod} we give a detailed description
of the moduli space $\wtM$ of thin sincere representations
via its ``defining polytope'' $\kD(\wt)$. 
The easiest way to obtain $\kD(\wt)$ from the quiver is to imagine
$Q$ as a one-way pipe system carrying liquid; a weight $\wt\in\Hh$ describes
the input (possibly negative) into the system at each knot. 
Using this language, $\kD(\wt)$ 
is simply given as the set of all possible flows 
respecting both the direction of 
the pipes and the given input $\wt$ (see \zitat{Mod}{3}).\\
Considering the opposite viewpoint, {\em each} flow through 
our pipe system requires 
a certain input, i.e.\ a weight. In particular, from the special flow that is
constant $1$ at each pipe we obtain a special, so-called 
canonical, weight $\cwt$. The corresponding $\kD(\cwt)$ is a reflexive
polytope (in the sense of Batyrev, \cite{Ba}), i.e.\ the moduli space
$\cwtM$ is Fano (Proposition \zitat{Mod}{7}). \\
Fixing a weight $\wt$ in general position, i.e.\ $\wtM$ is smooth, flows
and weights have still another 
meaning.  Each flow defines an equivariant, with respect to the
defining torus, effective divisor, and each weight $\wt'$ defines an element
$\cL(\wt')$ in the Picard group of $\wtM$. 
Assigning a flow its input weight
corresponds to assigning a divisor its class in the Picard group
(see \zitat{Uni}{1}). 
\par

{\bf Example:} 1) In the special case $\wt' := \wt$
this recovers our ample line bundle introduced before. \\
2) The line bundle $\UB_p^{-1} \otimes \UB_q$ corresponds
to the weight with values $1$ at $p$, $-1$ at $q$, and zero at all other
points. 
\par


\neu{Int-3}
Our first main result is Theorem \zitat{Uni}{6} stating the lack of
self-extensions of  
$\UB$ on the moduli space $\cwtM$ with respect to the canonical weight, i.e.\
$\gExt^{l}_{\cwtM}(\UB,\UB) = 0$ for all $l >0$. This is 
proved by using a slightly generalized Kodaira 
vanishing argument which works for toric varieties, 
cf.\ Theorem \zitat{Uni}{5}.   
As a Corollary of Theorem \zitat{Uni}{6} we conclude that we
obtain a full and faithful functor from the bounded derived category of
finitely generated right modules over the endomorphism algebra $\Alg$ of $\UB$
into the bounded derived category of coherent sheaves on the moduli space
$\cwtM$ (Theorem \zitat{End}{4}). 
Moreover, in Theorem \zitat{End}{3} we provide a criterion for $\Alg =
\gEnd_{\cwtM}(\UB,\UB)$ to be  isomorphic to the path algebra $\ck Q$ of the 
quiver $Q$. \\
Combining both results we obtain
the following relation between the derived categories of right $\ck Q$-modules
and of coherent sheaves on $\cwtM$, respectively{\DP}
\par

{\bf Theorem:}
{\em Assume $Q$ is a quiver lacking $(1,0)$- and $(t,t)$-walls (see 
\zitat{Mod}{2} for an explanation). 
Then, 
\[
-\otimes_{\ck Q}^{\Ll}\UB: \cD^b\Big(\mmodstr\ck Q \Big) \lra \cD^b \Big(
\mCoh\big( \cwtM \big)\Big)
\]
is a full and faithful functor from the bounded derived category
of finitely generated right $\ck Q$-modules into the bounded derived category
of coherent sheaves on $\cwtM$.
}
\par


\neu{Int-4}
The result above is closely related to tilting theory. Since the
fundamental paper \cite{Beilinson}, tilting theory has become a major tool in
classifying vector bundles; a tilting sheaf induces an equivalence of
bounded derived categories, as in the previous Theorem. 
To be precise we recall the 
definition of a tilting sheaf (\cite{Baer}). A sheaf $\cT$ on a smooth
projective variety is called a {\sl tilting sheaf} if
\kitem
\vspace{-1ex}
\item[1)] it has no higher self-extensions, that is $\mExt^l(\cT,\cT)=0$ for all
$l >0$,
\vspace{-1.5ex} 
\item[2)] the direct summands generate the bounded derived category, and
\vspace{-1.5ex}
\item[3)] the endomorphism algebra $\Alg$ of $\cT$ has finite global (homological)
dimension. 
\vspace{-1.5ex}
\kenditem
Then, the functors $\Rr \mHom(\cT,-)$ and
$-\otimes^{\Ll}_{\Alg}\cT$ define mutually inverse equivalences of the
bounded derived categories of coherent sheaves on the underlying variety of
$\cT$ and of the finitely generated right $\Alg$-modules, respectively. 
For constructions of tilting
bundles and their relations to derived categories we refer to the following
papers \cite{Kapranov}, \cite{Beilinson}, \cite{Rudakov}, \cite{Bondal},
and \cite{Orlov}. For the similar notion of a tilting module we refer to
\cite{HRTilt}.\\ 
For our purpose, the notion of an exceptional sequence is more useful. Let
$\cC$ be any of the categories introduced above: the category of finitely
generated right modules over a finite dimensional algebra, the category of
coherent sheaves on a smooth projective variety, or one of its derived
categories. Thus, $\cC$ is either 
an abelian or a triangulated $\ck$-category. Each object in $\cC$ has a
unique, up to isomorphism and reordering, decomposition into indecomposable
direct summands, i.e.\ $\cC$ is a  Krull-Schmidt category. Moreover, the
extension groups are defined and globally bounded; they are finite-dimensional
$\ck$-vector spaces. 
An object in $\cC$ is called {\sl exceptional} if it
has no self-extensions and its endomorphism ring is $\ck$.
A sequence $(\cE_0,\ldots,\cE_n)$ of
objects in $\cC$ is called {\em exceptional} if 
\kitem
\vspace{-1ex}
\item[1)] each object $\cE_i$ for $i=0,\ldots,n$ is exceptional and 
\vspace{-1.5ex}
\item[2)] $\mExt^l(\cE_j,\cE_i) = 0$ for all $l\geq
0$, and $j>i$. 
\vspace{-1.5ex}
\kenditem
Such a sequence is called {\sl strong exceptional} if,
additionally,
\kitem
\vspace{-1ex}
\item[3)] $\mExt^l(\cE_i,\cE_j) = 0$ for all $l>0$ and all $i,j
=0,\ldots,n$. 
\vspace{-1.5ex}
\kenditem
Finally, it is called {\sl full} if in addition to 1), 2) and 3)
\kitem
\vspace{-1ex}
\item[4)] the objects
$\cE_i$ for $i=0,\ldots,n$ generate the bounded derived category. 
\vspace{-1.5ex}
\kenditem
Thus, each
full strong exceptional sequence defines a tilting bundle $\oplus_{i=0}^n
\cE_i$, because the endomorphism algebra of $\oplus_{i=0}^n \cE_i$ has
global dimension at most $n$. 
Vice versa, each tilting bundle whose direct summands are line bundles 
gives rise to a strong exceptional sequence.\\
Using this language, our vanishing result Theorem \zitat{Uni}{6} means
that the direct summands of $\UB$ form a strong exceptional sequence. 
\par


\neu{Int-5}
In general, this sequence cannot be full. Assume the contrary; then the
bounded derived categories in the previous theorem are equivalent. The
first one is a derived category of a hereditary abelian category, whose
structure is well-known (\cite{RinHer}). In particular, the Serre functor
(see \cite{BK} for the definition)
coincides with the Auslander-Reiten translation and fixes objects up to
translation only in case the category is tame or just semi-simple
(cf.\ \cite{Happel} \S 1.4/5). On the
other hand, the Serre functor in the bounded derived 
category of coherent sheaves fixes all skyscraper sheaves up to a
translation. Consequently, an  equivalence implies that $\cM$ is a point or
a projective line in case the algebra $\ck Q$ is semi-simple or tame,
respectively.  It follows that $Q$
is a point or the Kronecker quiver; the remaining tame cases may not
appear (see \cite{Ringel} Theorem p.\ 158). \\
Nevertheless, there is some hope that
one may find a complement $\overline{\UB}$ such that $\UB \oplus
\overline{\UB}$ is a tilting bundle. At least a class of very particular
examples of tilting bundles on toric quiver varieties is known
(\cite{HilleTor}, Theorem 3.9). 
\par



\neu{Int-6} 
For an introduction to quivers and path algebras we refer 
the reader to \cite{Ringel}
and \cite{ARS}; the theory of localizations may be found in
\cite{Schofield}. For an introduction to moduli spaces we mention
\cite{Newstead} and for moduli of representations of quivers we refer to
the work of King \cite{King}. For results on triangulated categories we
refer to \cite{Happel} and \cite{HartshorneRes}. Our standard reference for
toric geometry is \cite{Ke}; for a short introduction to this area we
also mention \cite{Fulton}.\\
We would like to thank G.~Hein, A.~King, and A.~Schofield for helpful 
discussions.
\par

%
%

\section{Moduli spaces of thin sincere representations}\label{Mod}


\neu{Mod-1}
Let $Q$ be a connected quiver without oriented cycles; it consists of a set
$Q_0$ of vertices, a set $Q_1$ of arrows, and two functions $s,t:Q_1 \ra
Q_0$ assigning to each arrow $\alpha \in Q_1$ its source $s(\alpha)$ and
its target $t(\alpha)$. A representation of $Q$ is a collection of finite
dimensional $\ck$-vector spaces $x(q)$ for each vertex $q$ 
together with a collection
of linear maps $x(\alpha): x(s(\alpha)) \ra x(t(\alpha))$ for each arrow
$\alpha\in Q_1$. The dimension vector $d = (d_q \mid q \in Q_0)$ of a
representation $x$ is defined by $d_q = \dim x(q)$. A representation is
called {\sl thin} if $\dim x(q) \leq  1$ for all $q \in Q_0$ and {\sl
sincere} if $\dim x(q) \geq 1$ for all $q \in Q_0$. In this paper we
consider only thin sincere representations. \\
We denote by $\cR = \oplus_{\alpha \in Q_1} \ck$ the space of all thin
sincere representations, that is $x(q) = \ck$. By 
$G = \times_{q\in Q_0} \ck^{\ast}$ we denote the torus
acting via conjugation on $\cR$. The orbits of this action are exactly the
isomorphism classes of thin sincere representations, i.e.\ their moduli space
may be obtained via GIT. Doing so, we have to deal with the notion of
stability with respect to a given weight (cf. \cite{King}){\DP}
\par

{\bf Definition:}
The elements of the real vector space 
$\Hh := \{\wt: Q_0 \ra \Rr\mid \sum_{q \in Q_0} \wt_q = 0 \}$ 
are called {\sl weights} of the quiver $Q$.\\
Let $\wt\in\Hh$. 
A thin sincere representation $x$ of $Q$ is {\sl $\wt$-stable} ({\sl
-semistable}) if for
each proper non-trivial subrepresentation $y \subset x$ we have 
$\sum_{q \in Q_0\mid y(q) \not= 0} \wt_q <0$ ($\leq 0$ respectively). Two
semistable representations $x$ and $y$ are called {\sl $S$-equivalent} with
respect to $\wt$ if 
the factors of the stable Jordan-H\"older filtration coincide.\\
A subquiver $Q' \subseteq Q$ with $Q'_0=Q_0$ is
{\sl $\wt$-stable} ({\sl $\wt$-semistable}) if it has a $\wt$-stable
($\wt$-semistable) representation. Two quivers are {\sl $S$-equivalent} with
respect to $\wt$ if they admit $\wt$-semistable representations of the
same $S$-equivalence classes.\\  
Finally, we denote by $\cT(\wt)$ the set of all $\wt$-semistable subtrees $T
\subseteq Q$ with $T_0 = Q_0$ and by
$Q_1(\wt)$ the set of all arrows $\alpha$ such that 
$Q \setminus \{\alpha\}$ is a $\wt$-stable subquiver.
\par

In other words, a representation $x$ is $\wt$-stable precisely when  the
subquiver  $Q_0 \cup \{ \alpha \in Q_1 \mid  x(\alpha) \not= 0 \}$ is
$\wt$-stable. Moreover, a subquiver $Q'$ is $\wt$-stable if and only if for
all non-trivial proper subsets $S \subset Q_0$ which are closed under
successors in $Q'$ we have $\sum_{q \in S} \wt_q < 0$ (see \cite{HilleThin}
\S 2 and \cite{HilleTor} Lemma 1.4). We also note that each $S$-equivalence
class contains a 
unique minimal $\wt$-semistable subquiver -- just take the disjoint union of
the support of the Jordan-H\"older factors. 
\par


\neu{Mod-2}
As already mentioned in the beginning, for any given weight $\wt$ the moduli
space $\wtM$ exists; however, different weights $\wt$ may cause different 
moduli spaces. According to \cite{HilleThin} there is a chamber system in
$\Hh$, and the type of $\wtM$ can only flip if $\wt$ crosses walls of the
following type:
\par

{\bf Definition:}
$W\subseteq\Hh$ is called a $(t^+,t^-)${\sl -wall} if 
\[
W=\Big\{ \wt \in \Hh \mid \sum_{q \in Q_0^+} \wt_q= -\sum_{q \in Q_0^-}
\wt_q=0 \Big\} 
\]
for some decomposition 
$Q_0 = Q_0^+ \sqcup Q_0^-$ such that the full subquivers $Q^+$ and $Q^-$
are both connected and such there are exactly $t^+$ arrows pointing from
$Q_0^+$ to $Q_0^-$ and $t^-$ arrows the other way around.
We say that $\wt$ is in {\em general position} if $\wt$ does not lie on any wall 
and if the moduli space is not empty. 

Assume $\wt$ is in general position, then $\wtM$ is smooth and has the
(maximal) dimension $d = \#Q_1\, -\, \#Q_0\, + \,1$, see
\cite{HilleThin}. Moreover, for those weights, every semistable thin
sincere representation is stable. 
\par


\neu{Mod-3}
To describe the toric structure of $\wtM$ we introduce the 
real vector space of flows defined as $\Fh:= \{r: Q_1\to \R\}=\R^{Q_1}$.
A flow is called {\sl regular} if it has only non-negative values, i.e.\ if it
respects the direction of the pipes. For any $\alpha\in Q_1$ we denote by
$f^\alpha\in\Fh$ the characteristic flow mapping $\alpha$ to $1$ and keeping
the remaining pipes dry. More generally, for each walk $w$ 
without cycles in $Q$ we define
the characteristic flow $f^w$ mapping an arrow $\alpha \in w$ to $1$,
an arrow $\beta$ with $\beta^{-1} \in w$ to $-1$, and the remaining arrows
to $0$. This flow is regular if and only if $w$ is a path, i.e.\ respects
the orientation in $Q$. \\
There is a canonical linear map $\pi:\Fh\to\Hh$ describing the input of flows;
if $\alpha\in Q_1$ points from $p$ to $q$, then $\pi(f^\alpha)$ sends $p$ and
$q$ onto $1$ and $-1$, respectively. Thus
$$
(\pi(r))_q = \sum_{s(\alpha) = q} r_{\alpha} - \sum_{t(\alpha) = q}
r_{\alpha}. 
$$
This leads to the following definition:
\par

{\bf Definition:}
The convex {\em polytope of flows} $\kD(\wt)$ assigned to a weight $\wt$ 
is defined as the intersection
\[
\kD(\wt):=\pi^{-1}(\wt)\cap \R^{\kQ_1}_{\geq 0}\,.
\]
This means that $\kD(\wt)$ consists of exactly those regular flows respecting the
prescribed input $\wt$. Moreover, $\kD(\wt)$ is compact since $Q$ has no
oriented cycles.\\
The vector spaces $\Fh$ and $\Hh$ contain the lattices $\Fh_\Z$ and $\Hh_\Z$ 
of integral flows and weights, respectively. For any integral weight
$\wt:\kQ_0\to\Z$ we define the affine lattice
$\kMwt\subseteq\Fh$ as the fiber $\Fh_\Z\cap\pi^{-1}(\wt)$, i.e.\
\[
\kMwt = \Big\{r\in\Z^{\kQ_1}\,|\;
\sum\limits_{s(\alpha)=q}r_\alpha - \sum\limits_{t(\alpha)=q}r_\alpha
=\wt_q\;\mbox{ for all }q\in\kQ_0\Big\}\,.
\]
Any element $r \in \kMwt$ provides obviously an isomorphism
$(+r): M:=M^0 \stackrel{\sim}{\longrightarrow} \kMwt$.
\par


\neu{Mod-4}
The following lemma will be crucial for the understanding of our flow
polytope as well as for proving the upcoming vanishing theorem in \S \ref{Uni};
it explicitly provides points of the lattices $\kMwt${\DP}
Let $T\subseteq \kQ$ be an arbitrary
maximal tree. Each arrow
$\alpha\in T_1$ divides $\kQ_0$ into two disjoint 
subsets, the source $S_T(\alpha)$ and the target $T_T(\alpha)$.
\par

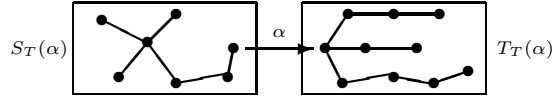
\begin{center}
\setlength{\unitlength}{0.00600in}%
\begingroup\makeatletter
\def\x#1#2#3#4#5#6#7\relax{\def\x{#1#2#3#4#5#6}}%
\expandafter\x\fmtname xxxxxx\relax \def\y{splain}%
\ifx\x\y   
\gdef\SetFigFont#1#2#3{%
  \ifnum #1<17\tiny\else \ifnum #1<20\small\else
  \ifnum #1<24\normalsize\else \ifnum #1<29\large\else
  \ifnum #1<34\Large\else \ifnum #1<41\LARGE\else
     \huge\fi\fi\fi\fi\fi\fi
  \csname #3\endcsname}%
\else
\gdef\SetFigFont#1#2#3{\begingroup
  \count@#1\relax \ifnum 25<\count@\count@25\fi
  \def\x{\endgroup\@setsize\SetFigFont{#2pt}}%
  \expandafter\x
    \csname \romannumeral\the\count@ pt\expandafter\endcsname
    \csname @\romannumeral\the\count@ pt\endcsname
  \csname #3\endcsname}%
\fi
\endgroup
\begin{picture}(395,80)(55,680)
\thicklines
\put(300,720){\circle*{10}}
\put(335,720){\circle*{10}}
\put(400,750){\circle*{10}}
\put(360,695){\circle*{10}}
\put(320,750){\circle*{10}}
\put(425,700){\circle*{10}}
\put(395,690){\circle*{10}}
\put(315,690){\circle*{10}}
\put(360,750){\circle*{10}}
\put(380,720){\circle*{10}}
\put(105,745){\circle*{10}}
\put(170,750){\circle*{10}}
\put(170,690){\circle*{10}}
\put(120,695){\circle*{10}}
\put(145,725){\circle*{10}}
\put(215,695){\circle*{10}}
\put(220,720){\circle*{10}}
\thinlines
\put( 80,760){\line( 0,-1){ 80}}
\put( 80,680){\line( 1, 0){160}}
\put(240,680){\line( 0, 1){ 80}}
\put(240,760){\line(-1, 0){160}}
\put(280,760){\line( 0,-1){ 80}}
\put(280,680){\line( 1, 0){160}}
\put(440,680){\line( 0, 1){ 80}}
\put(440,760){\line(-1, 0){160}}
\thicklines
\put(450,715){\makebox(0,0)[lb]{\smash{\SetFigFont{12}{14.4}{rm}\scriptsize $T_T(\alpha)$}}} 
\put(230,720){\vector( 1, 0){ 60}}
\put(260,730){\makebox(0,0)[cc]{\smash{\SetFigFont{12}{14.4}{rm}\scriptsize $\alpha$}}}
\put(105,745){\line( 2,-1){ 40}}
\put(145,725){\line( 1, 1){ 25}}
\put(145,725){\line( 3,-4){ 25.800}}
\put(170,690){\line( 6, 1){ 44.595}}
\put(215,695){\line( 1, 5){  5}}
\put(120,695){\line( 5, 6){ 25}}
\put(300,720){\line( 1, 0){ 35}}
\put(335,720){\line( 1, 0){ 45}}
\put(300,720){\line( 2, 3){ 20}}
\put(320,750){\line( 1, 0){ 40}}
\put(360,750){\line( 1, 0){ 40}}
\put(300,720){\line( 1,-2){ 15}}
\put(315,690){\line( 6, 1){ 44.595}}
\put(360,695){\line( 6,-1){ 34.865}}
\put(395,690){\line( 3, 1){ 30}}
\put( 25,715){\makebox(0,0)[lb]{\smash{\SetFigFont{12}{14.4}{rm}\scriptsize $S_T(\alpha)$}}}
\end{picture}
\end{center}
\par

{\bf Lemma:}
{\em
Fix a tree $T$ and let $\wt\in\Hh$. For any flow $\keps$ there is a unique
element $r = r^T \in \pi^{-1}(\wt) \subseteq \Fh$ satisfying
$r_\alpha := \keps(\alpha)$ for $\alpha \notin T_1$. Its remaining coordinates
(i.e.\ for $\alpha\in T_1$) are given by
\[
r_\alpha\,-\,\keps(\alpha)\;=\;
\sum_{q\in S_T(\alpha)}\!\!\wt_q\;
+ \sum_{T_T(\alpha)\stackrel{\beta}{\rightarrow}S_T(\alpha)}
\!\!\!\!\keps(\beta)\;
- \sum_{S_T(\alpha)\stackrel{\beta}{\rightarrow}T_T(\alpha)}
\!\!\!\!\keps(\beta)\,.
\vspace{-2ex}
\]
}
\par

{\bf Proof:}
First, we should note that both $\keps$ and $r$ are flows -- the different
notation for their coordinates ($\keps(\alpha)$ and $r_\alpha$, respectively)
was chosen for psychological reasons only.
Now, let $r$ be {\em some} element of $\Fh$. If $\pi(r)=\wt$, then for 
any subdivision $Q_0 = Q_0^+ \sqcup Q_0^-$ we obtain
\[
\sum_{s(\beta) \in Q_0^+, t(\beta) \in Q_0^-}  r_\beta -
\sum_{s(\beta) \in Q_0^-, t(\beta) \in Q_0^+}  r_\beta = \sum_{q
\in Q_0^+}\wt_q = -\sum_{q \in Q_0^-}\wt_q
\]
just by summing up the $\kMwt$-equations with $q\in Q_0^+$. The reverse
implication is also true, even if we restrict ourselves to the special
subdivisions provided by arrows $\alpha\in T_1$ via 
$Q_0^+:= S_T(\alpha)$ and $Q_1^+:= T_T(\alpha)$. \\
On the other hand, these subdivisions have the important property that
$\alpha$ is the only arrow that belongs to both the tree $T$ and to one of the
index sets $\{\beta\,|\; s(\beta) \in Q_0^+, t(\beta) \in Q_0^-\}$ or
$\{\beta\,|\; s(\beta) \in Q_0^-, t(\beta) \in Q_0^+\}$. In particular, 
by just taking care of this single exception, in
the above equations we may always replace $r_\beta$ by $\keps(\beta)$. 
\hfill$\Box$
\par

If the weight $\wt$ and the flow 
are integral, then so is $r$, i.e.\ $r\in \kMwt$.
\par


\neu{Mod-5}
{\bf Proposition:}
{\em
Let $\wt\in\Hh$ be an integral weight, then
the polytope of flows $\kD(\wt)\subseteq\kMwt$ is always a lattice polytope.
The associated projective toric variety equals $\wtM$. Moreover, $\kD(\wt)$
provides an ample, equivariant line bundle $\cL(\wt)$ on $\wtM$.}
\par

{\bf Proof:}
First, we establish a
one-to-one correspondence between vertices of $\kD(\wt)$ and $S$-equivalence
classes of $\wt$-semistable trees.
Faces of $\kD(\wt)$ in any dimension are obtained by forcing certain 
coordinates of $\Fh$ to be zero.
In particular, vertices are points with a
maximal set of vanishing coordinates. 
Let $r\in\kD(\wt)$ be a vertex and denote its support by
\[
\supp\, r:=\left\{\alpha\in Q_1\,| \; r_\alpha >0\right\}\,.
\] 
If $\supp\, r$ contained any cycle of $Q$, then we could replace $r$ by a
different regular flow with the same weight and a smaller support. Hence,
$\supp\, r$ is contained in maximal trees of $Q$. Moreover, we obtain that
\kitem
\item
every maximal tree $T$ containing $\supp\, r$ is $\wt$-semistable, and
\item
those trees are stable if and only if $\supp\, r=T_1$ (which
determines the tree uniquely).
\kenditem
To prove these facts, take a proper subset $Q_0^+\subset Q_0$ that is closed
under successors in $T$; denoting $Q_0^-:=Q_0\setminus Q_0^+$ this means
that there are no arrows pointing from $Q_0^+$ to $Q_0^-$ in $T$. Hence,
using the formula mentioned in the previous proof, 
$$
\sum_{q \in Q_0^+}\wt_q = - \sum_{[Q_0^-\stackrel{\beta}{\to} Q_0^+]}
r_\beta \leq 0. 
$$

Conversely, let $T$ be any maximal
tree. The previous lemma tells us that there is exactly one
$r\in\pi^{-1}(\wt)\subseteq \Fh$ 
such that $\supp\, r\subseteq T_1$; it has integer coordinates 
$r_\alpha=\sum_{S_T(\alpha)}\wt_q$ for $\alpha\in T_1$. Moreover, 
if $T$ is $\wt$-semistable, then these numbers are non-negative, meaning
that $r\in\kD(\wt)$. Thus, it must be a vertex. Moreover, two different
trees define the same vertex if and only if they are $S$-equivalent.
\par

What does the inner normal fan $\kSwt$ look like?
Denote by $\kA$ the matrix describing the incidences of our quiver;
$\kA$ consists of $\# \kQ_0 $ rows and $\# \kQ_1 $ columns, and for
$q\in \kQ_0$, $\alpha\in\kQ_1$ we have
\[
\kA_{q\alpha}:=\left\{
\begin{array}{cl}
+1&\mbox{if $q = s(\alpha)$ }\\
-1&\mbox{if $q = t(\alpha)$ }\\
0&\mbox{otherwise}\,.
\end{array}
\right.
\]
The free abelian group
$N:=\Z^{\kQ_1}/ (\kA\mbox{--rows})$ is dual to the lattice 
$M:=M^0=\big(\kA\mbox{--rows}\big)^\bot\subseteq\Fh$. Hence,
writing $a^\alpha\in N$ for the images of the canonical
vectors $e^\alpha\in\Z^{\kQ_1}$, we obtain
$$
\kSwt=\Big\{ \langle a^{\alpha_1},\dots,a^{\alpha_k}\rangle \, \Big|
\begin{array}[t]{l}
Q \setminus \{\alpha_1,\dots,\alpha_k\} \mbox{ is the minimal
element in some $S$-equivalence class }\\
\mbox{ of $\wt$-semistable subquivers } \Big\}.
\end{array}
$$
According to \cite{HilleTor}, Theorem 1.7 and (3.3), this is exactly the fan 
defining the moduli space $\wtM$.  
\hspace*{\fill}$\Box$
\par

{\bf Example:} In the quiver below the canonical weight
is $\cwt = (2,1,-1,-2)$. 
$$
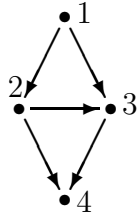
\setlength{\unitlength}{0.00600in}%
\begingroup\makeatletter
\def\x#1#2#3#4#5#6#7\relax{\def\x{#1#2#3#4#5#6}}%
\expandafter\x\fmtname xxxxxx\relax \def\y{splain}%
\ifx\x\y   
\gdef\SetFigFont#1#2#3{%
  \ifnum #1<17\tiny\else \ifnum #1<20\small\else
  \ifnum #1<24\normalsize\else \ifnum #1<29\large\else
  \ifnum #1<34\Large\else \ifnum #1<41\LARGE\else
     \huge\fi\fi\fi\fi\fi\fi
  \csname #3\endcsname}%
\else
\gdef\SetFigFont#1#2#3{\begingroup
  \count@#1\relax \ifnum 25<\count@\count@25\fi
  \def\x{\endgroup\@setsize\SetFigFont{#2pt}}%
  \expandafter\x
    \csname \romannumeral\the\count@ pt\expandafter\endcsname
    \csname @\romannumeral\the\count@ pt\endcsname
  \csname #3\endcsname}%
\fi
\endgroup
\begin{picture}(100,175)(150,590)
\thicklines
\put(160,680){\circle*{10}}
\put(240,680){\circle*{10}}
\put(200,600){\circle*{10}}
\put(205,750){\vector( 1,-2){ 30}}
\put(195,750){\vector(-1,-2){ 30}}
\put(200,760){\circle*{10}}
\put(170,680){\vector( 1, 0){ 60}}
\put(210,590){\makebox(0,0)[lb]{\smash{\SetFigFont{12}{14.4}{rm}4}}}
\put(165,670){\vector( 1,-2){ 30}}
\put(235,670){\vector(-1,-2){ 30}}
\put(210,755){\makebox(0,0)[lb]{\smash{\SetFigFont{12}{14.4}{rm}1}}}
\put(150,690){\makebox(0,0)[lb]{\smash{\SetFigFont{12}{14.4}{rm}2}}}
\put(250,675){\makebox(0,0)[lb]{\smash{\SetFigFont{12}{14.4}{rm}3}}}
\end{picture}
 \hspace{3cm}
\begin{array}{c}
C = \left(
\begin{array}{rrrrr}
1 & 1 & 0 & 0 & 0 \\
-1& 0 & 1 & 1 & 0 \\
0 & -1& -1& 0 & 1 \\
0 & 0 & 0 & -1& -1
\end{array}
\right) \\ ~ \\ ~ \\ ~ \\ ~ \\ ~ \\ ~ 
\end{array}
$$

\vspace*{-2cm} 

The corresponding fan $\kSwt$, and the polytope $\Delta(\cwt)$ look 
like the following:

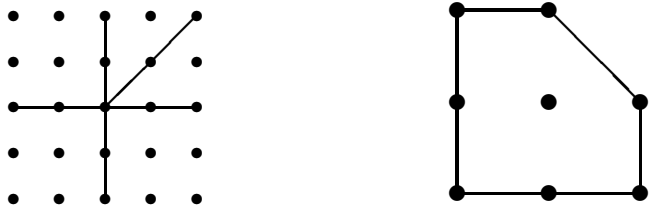
\begin{center}
\setlength{\unitlength}{0.00600in}%
\begingroup\makeatletter
\def\x#1#2#3#4#5#6#7\relax{\def\x{#1#2#3#4#5#6}}%
\expandafter\x\fmtname xxxxxx\relax \def\y{splain}%
\ifx\x\y   
\gdef\SetFigFont#1#2#3{%
  \ifnum #1<17\tiny\else \ifnum #1<20\small\else
  \ifnum #1<24\normalsize\else \ifnum #1<29\large\else
  \ifnum #1<34\Large\else \ifnum #1<41\LARGE\else
     \huge\fi\fi\fi\fi\fi\fi
  \csname #3\endcsname}%
\else
\gdef\SetFigFont#1#2#3{\begingroup
  \count@#1\relax \ifnum 25<\count@\count@25\fi
  \def\x{\endgroup\@setsize\SetFigFont{#2pt}}%
  \expandafter\x
    \csname \romannumeral\the\count@ pt\expandafter\endcsname
    \csname @\romannumeral\the\count@ pt\endcsname
  \csname #3\endcsname}%
\fi
\endgroup
\begin{picture}(170,170)(235,475)
\thicklines
\put(280,560){\circle*{10}}
\put(240,560){\circle*{10}}
\put(320,520){\circle*{10}}
\put(280,520){\circle*{10}}
\put(240,520){\circle*{10}}
\put(240,480){\circle*{10}}
\put(280,480){\circle*{10}}
\put(320,480){\circle*{10}}
\put(360,480){\circle*{10}}
\put(360,520){\circle*{10}}
\put(360,560){\circle*{10}}
\put(320,560){\circle*{10}}
\put(360,600){\circle*{10}}
\put(280,600){\circle*{10}}
\put(320,600){\circle*{10}}
\put(240,600){\circle*{10}}
\put(400,640){\line(-1,-1){ 80}}
\put(240,640){\circle*{10}}
\put(280,640){\circle*{10}}
\put(320,640){\circle*{10}}
\put(360,640){\circle*{10}}
\put(400,640){\circle*{10}}
\put(400,600){\circle*{10}}
\put(400,560){\circle*{10}}
\put(400,520){\circle*{10}}
\put(400,480){\circle*{10}}
\put(320,560){\line(-1, 0){ 80}}
\put(320,560){\line( 0,-1){ 80}}
\put(320,560){\line( 1, 0){ 80}}
\put(320,560){\line( 0, 1){ 80}}
\end{picture}
\hspace{3cm}
\setlength{\unitlength}{0.01200in}%
\begingroup\makeatletter
\def\x#1#2#3#4#5#6#7\relax{\def\x{#1#2#3#4#5#6}}%
\expandafter\x\fmtname xxxxxx\relax \def\y{splain}%
\ifx\x\y   
\gdef\SetFigFont#1#2#3{%
  \ifnum #1<17\tiny\else \ifnum #1<20\small\else
  \ifnum #1<24\normalsize\else \ifnum #1<29\large\else
  \ifnum #1<34\Large\else \ifnum #1<41\LARGE\else
     \huge\fi\fi\fi\fi\fi\fi
  \csname #3\endcsname}%
\else
\gdef\SetFigFont#1#2#3{\begingroup
  \count@#1\relax \ifnum 25<\count@\count@25\fi
  \def\x{\endgroup\@setsize\SetFigFont{#2pt}}%
  \expandafter\x
    \csname \romannumeral\the\count@ pt\expandafter\endcsname
    \csname @\romannumeral\the\count@ pt\endcsname
  \csname #3\endcsname}%
\fi
\endgroup
\begin{picture}(90,90)(275,515)
\thicklines
\put(280,560){\circle*{7}}
\put(320,520){\circle*{7}}
\put(280,520){\circle*{7}}
\put(360,520){\circle*{7}}
\put(360,560){\circle*{7}}
\put(320,600){\circle*{7}}
\put(320,560){\circle*{7}}
\put(280,520){\line( 0, 1){ 80}}
\put(280,600){\circle*{7}}
\put(280,600){\line( 1, 0){ 40}}
\put(320,600){\line( 1,-1){ 40}}
\put(360,560){\line( 0,-1){ 40}}
\put(360,520){\line(-1, 0){ 80}}
\end{picture}
\end{center}

It is known that the toric variety of this fan is the blow
up of the projective plane in two points, which is isomorphic to the blow up
of the two-dimensional smooth quadric in one point. 
\par


\neu{Mod-6}
Equivariant (with respect to the torus action), invertible sheaves ${\cL}$ on a
toric variety $X(\Sigma)$ are completely determined by their order function
$\mbox{ord}\,{\cL}:\Sigma^{(1)}\to \Z$ or its piecewise linear continuation
$\mbox{ord}\,{\cL}:N_{\R}\to\R$. If $r(\sigma) \in M$ provides a local
generator $x^{r(\sigma)}$
of ${\cL}$ on $U(\sigma) \subseteq X(\Sigma)$, then $\mbox{ord}\,{\cal L}(a)$
is defined as $\langle a, r(\sigma) \rangle$ if $a\in\sigma$ (cf.\ \cite{Oda}).
Moreover, if ${\cal L}$ is an ample (or at least globally generated)
invertible sheaf given by a lattice polytope  
$\kD\subseteq M_{\R}$, then the local generators of ${\cL}$ correspond
to the  vertices of $\kD$. In particular, 
$\mbox{ord}\,{\cL}(a)=\mbox{min}\,\langle a, \kD\rangle$. Shifting the
polytope 
$\kD$ by a vector $r\in M$ means to replace ${\cL}$ by $x^{r}\cdot
{\cal L}$.  
The corresponding
order functions differ by the globally linear function
$\langle \bullet ,r\rangle$.
\par

{\bf Lemma:}
{\em
Let $r\in\kMwt$ be an arbitrary element. Then, the mapping
$a^\alpha\mapsto -r_\alpha$ gives the order function of $\cL(\wt)$ on
$\wtM$. A different choice $r^\prime\in\kMwt$ just changes the order function
by the linear summand $\langle \bullet, r' - r \rangle$.}
\par

{\bf Proof:}
We may use $r\in\kMwt$ to carry $\kD(\wt)$ into the ``right'' lattice 
$M$ (see the end of \zitat{Mod}{3}). Then,  
the order function applied to $a^\alpha\in\kSwt^{(1)}\subseteq N_{\R}$ is
\[
\mbox{ord}\,\cL(\wt)(a^\alpha)=\mbox{min}\,\langle
a^\alpha,\kD(\wt)-r\rangle 
= \mbox{min}\,\langle e^\alpha,\kD(\wt)\rangle -r_\alpha =-r_\alpha\,.
\vspace{-4ex}
\]
\hspace*{\fill}$\Box$
\par


\neu{Mod-7}
Given the quiver $\kQ$, the canonical weight $\cwt$ announced in the
introduction is defined as the weight of the flow $r^c$ that 
is constant $1$ on every arrow. Explicitly, this means
\[
\cwt(q):=\#\{\mbox{arrows with source } q\} - 
\#\{\mbox{arrows with target } q\}\,.
\]
The advantage of $\cwt$ is the existence of a unique interior lattice point
in the polytope
$\kD(\cwt)$: it is again the flow $r^c=[1,\dots,1]$. 
\par

{\bf Proposition:}
{\em
The polytope $\kD(\cwt)$ is reflexive (in the sense of \cite{Ba}), its
order function is $-1$ on the generators $a^\alpha\in \kScwt^{(1)}$, and
the ample divisor $\cL(\cwt)$ is anti-canonical.}
\par

{\bf Proof:} The three claims are synonymous, i.e.\ we just have to look at
the order function of $\cL(\cwt)$. Applying Lemma \zitat{Mod}{6} on $r^c =
\ku{1}$ yields the result. 
\hfill$\Box$
\par

%
%

\section{The cohomology of the universal bundle}\label{Uni}


\neu{Uni-1} From 
now on we assume that $\wt$ is an integral weight in general position,
i.e.\ $\wtM$ is 
a smooth variety. To each integral flow we associate a divisor as follows:
\[
f^\alpha\in\Fh_\Z \mapsto D_\alpha:=\{ x \in \wtM \mid x_{\alpha} = 0\}=
\left\{ \begin{array}{ll}
\ko{\orb (\alpha)} & \mbox{if } a^\alpha\in\Sigma^{(1)}=Q_1(\wt)\\
\emptyset & \mbox{otherwise}
\end{array}\right.
\]
with $\ko{\orb (\alpha)}$ denoting the closed orbit corresponding to the
one-dimensional cone $\alpha$. One obtains surjective maps 
$\Fh_{\Zz} \surj\mbox{\rm Div}\,\wtM$ from the space of integral flows
onto the equivariant divisors and, as a consequence, 
$\Hh_{\Zz} \surj\mbox{\rm Pic}\,\wtM$, $\wt' \mapsto \cL(\wt')$ from the
integral weights to the 
Picard group (see also \cite{HilleTor}, Theorem 2.3). Applying the map
$\pi: \Fh_{\Zz} \surj \Hh_{\Zz}$ means to assign a divisor its class in the
Picard group. \\
Copying the definition of \zitat{Mod}{3}, every weight $\wt'$ gives rise to
\[
\overline{\kD}(\wt'):= \{r\in \pi^{-1}(\wt')\,|\;
   r_\alpha\geq 0 \mbox{ for } \alpha\in Q_1(\wt)\}\,.
\]
Even if $\cL(\wt')$ is not ample on $\wtM$, the polytope
$\overline{\kD}(\wt')$ may still be used to describe the global sections:
\par

{\bf Proposition:}
{\em
The lattice points of $\overline{\kD}(\wt')$ provide 
a basis of the global sections of $\cL(\wt')$. Moreover,
if $Q_1(\wt')\subseteq Q_1(\wt)$, then both polytopes
$\kD(\wt')$ and $\overline{\kD}(\wt')$ coincide.
}
\par

{\bf Proof:} Given $\wt'$, we choose a flow 
$r\in M^{\wt'}$ providing the order function of a divisor in the class
defined by $\wt'$. The corresponding polytope of global sections is
contained in $M = M^0$; via the isomorphism $(+r): M^0 \ra M^{\wt'}$ it is
mapped onto $\ko{\kD}(\wt')$.\\
If $Q_1(\wt')\subseteq Q_1(\wt)$, then $\ko{\kD}(\wt')$ sits between
$\kD(\wt')$ and 
$\{r\in \pi^{-1}(\wt')\,|\;
   r_\alpha\geq 0 \mbox{ for } \alpha\in Q_1(\wt')\}$.
On the other hand, Proposition \zitat{Mod}{5} implies that the latter two
polytopes
are equal; its proof shows quite directly that the inequalities parametrized
by $Q_1\setminus Q_1(\wt')$ are redundant for the definition of $\kD(\wt')$.
\hspace*{\fill}$\Box$
\par


\neu{Uni-2}
Since $\wt$ is in general position, there is a universal bundle $\UB$
on $\wtM$; it splits into a direct sum $\UB= \oplus_{q \in Q_0}\UB_q$ of line
bundles. The direct summands $\UB_{p,q}:=\UB_p^{-1}\otimes \UB_q$ 
of $\underline{\lEnd}(\UB)$ have the following shape:
Choose a walk from $p$ to $q$ along (possibly reversed) arrows
$\alpha_1^{\keps(1)},\dots,\alpha_m^{\keps(m)}$, i.e.\ 
$\alpha_1,\dots,\alpha_m\in \kQ_1$ and $\keps(i)\in\{\pm 1\}$.
Then, denoting by $\CO(\alpha):=\CO(D_{\alpha})$ the sheaf corresponding
to the prime divisor $D_{\alpha}$, 
\[
\UB_{p,q}=\UB_p^{-1}\otimes \UB_q=
\bigotimes\limits_{i=1}^m\,\CO(\alpha_i)^{\keps(i)} \,.
\]
In the Picard group of $\wtM$ this sheaf does not depend on
the particular choice of the walk from $p$ to $q$: Using the language of
\zitat{Uni}{1}, the sheaves $\otimes_i\,\CO(\alpha_i)^{\keps(i)}$ are induced 
from the flows $\sum_i\keps(i)\cdot f^{\alpha_i}$, 
which all have the same weight. 
\par

{\em Notation:}
Setting $\keps(\alpha^i):=\keps(i)$ and
$\keps(\alpha):=0$ for $\alpha\notin\{\alpha^1,\dots,\alpha^m\}$ provides
a function $\keps:\kQ_1\to\{1,-1,0\}$ for every walk. This is the
characteristic flow introduced in \zitat{Mod}{3}. Then, the sheaf
$\UB_{p,q}$ may be written as
$\UB_{p,q}=\UB(\keps)=\otimes_{\alpha\in\kQ_1}\CO(\alpha)^{\keps(\alpha)}$;
the corresponding weight $\wt_{p,q}:=\pi(\keps)$ has value $1$ in $p$,
$-1$ in $q$, and $0$ in all other vertices.
\par


\neu{Uni-3}
{\bf Proposition:}
{\em 
\kitem
\item[(1)]
Let $\wt\in\Hh_\Z$ be an integral weight in general position.
Then, the sheaves $\UB_{p,q}\otimes \cL(\wt)$ and
$\UB_{p,q}^{-1}\otimes \cL(\wt)$ are generated by their 
global sections.
\item[(2)]
If, additionally, $\wt=\cwt$, then the polytopes 
$\ko{\kD}(\cwt\pm \wt_{p,q})$ (describing the global sections) 
have the same dimension as $\kD(\cwt)$.
\kenditem}
\par

{\bf Proof:} 
Since $\,\UB_{q,p}=\UB_{p,q}^{-1}$, 
it is sufficient to consider the latter sheaf.
The corresponding polytope $\ko{\kD}(\wt-\wt_{p,q})$ may be studied in
different level sets:
\[
\begin{array}{rcl}
\ko{\kD}(\wt-\wt_{p,q}) &=&
\{r\in\pi^{-1}(\wt-\wt_{p,q})\,|\; r_\alpha\geq 0 \mbox{ for }
\alpha\in Q_1(\wt)\}\\
&\cong&
\{r\in\pi^{-1}(\wt)\,|\; r_\alpha\geq \keps(\alpha)\mbox{ for } 
\alpha\in Q_1(\wt)\}\,.
\end{array}
\]
We will use the second description. 
\vspace{1ex}\\
(1)
The vertices of $\kD(\wt)$ and thus also the top-dimensional cones of
$\kSwt$ are in a one-to-one correspondence with the
$\wt$-stable trees in $Q$.
Let $T\in\kTh$; the corresponding vertex $\kD^T$ of $\kD(\wt)$ provides a
local generator of $\cL(\wt)$. Since the $\kD^T$ are characterized by the 
property
$\kD^T_\alpha=0$ for $\alpha\notin T$, we obtain the local generators of
$\UB_{p,q}^{-1}\otimes \cL(\wt)$ 
from the lattice points $r^T\in\kMwt$ assigned via
Lemma \zitat{Mod}{4} to the map $\keps:\kQ_1\to\Z$
describing a walk from $p$ to $q$.\\
We have to show that these local generators $r^{T}$
are regular on any open, affine subset
$U_{T^\prime}\subseteq\wtM$ corresponding to some possibly different 
tree $T^\prime\in\kTh$.
That means,
it remains
to check that $r^T\in r^{T^\prime}+(\sigma_{T'})^{\veee}$ where
$\sigma_{T'}$ denotes the cone corresponding to $T'$, i.e.\
$\sigma_{T'}$ is spanned by those arrows {\em not} contained in $T^\prime$.
\vspace{1ex}\\
{\em Claim: Let $T\in\kTh$ be a $\wt$-stable tree, and let $\alpha$ be any
arrow in $\kQ$. Then $r^T_\alpha\geq \keps(\alpha)$ is true.}
\vspace{0.5ex}\\
Before we prove that claim, we remark that it solves our
problem, 
as, for any tree, we know for $T^\prime$ that
$r^{T^\prime}_\alpha=\keps(\alpha)$ for $\alpha\notin T^\prime$.
Hence, the claim implies
$r^T_\alpha\geq r^{T^\prime}_\alpha$ for $\alpha\notin T^\prime$. On
the other 
hand, $(\sigma_{T'})^{\veee}$ is just given by the inequalities
$r_\alpha\geq 0$ for those $\alpha$.
\vspace{1ex}\\
The claim is trivial for $\alpha\notin T$;
if $\alpha\in T$, we use the formula for $r_\alpha-\keps(\alpha)$ presented
in Lemma \zitat{Mod}{4}. First, as already used in \zitat{Mod}{5}, stability
of $T$ implies $\sum_{q\in S_T(\alpha)}\wt_q\geq 1$. Now, the point is to
interpret the two remaining sums well: together they just count the
number of  arrows $\alpha_i^{\keps(i)}$ in the walk from $p$ to $q$
pointing from $T_T(\alpha)$ to 
$S_T(\alpha)$ minus those from $S_T(\alpha)$ to $T_T(\alpha)$. In particular,
\[
\sum\limits_{\makebox[4em]{$\ks T_T(\alpha)\stackrel{\beta}{\rightarrow}
S_T(\alpha)$}} \keps(\beta)\;
- \sum\limits_{\makebox[4em]{$\ks S_T(\alpha)\stackrel{\beta}{\rightarrow}
T_T(\alpha)$}} \keps(\beta)
\;=\;
\left\{\begin{array}{rl}
-1 & \mbox{if } p\in S_T(\alpha) \mbox{ and } q\in T_T(\alpha)\\
1  & \mbox{if } q\in S_T(\alpha) \mbox{ and } p\in T_T(\alpha)\\
0  & \mbox{if } p,q\in S_T(\alpha) \mbox{ or } p,q\in T_T(\alpha)\,.
\end{array}\right.
\]
In any case, $r_\alpha-\keps(\alpha)$ remains non-negative.
\vspace{1ex}\\
(2)
For the second part, we consider $\wt=\cwt$. Since $\keps$ has only
$-1$, $0$, or $1$ as values, the canonical flow 
$r^c=\ku{1}\in\kD(\cwt)$ is also contained in 
$\ko{\kD}(\cwt-\wt_{p,q})\subseteq\pi^{-1}(\cwt)$.
Assuming $\dim \ko{\kD}(\cwt-\wt_{p,q}) < \dim \kD(\cwt)$, this means that
there exist arrows $\alpha^1,\dots,\alpha^k\in Q_1(\cwt)$ having the following
two properties:
\kitem
\item[(i)]
The flow $r^c$ satisfies the 
$\alpha^v$-inequalities of $\ko{\kD}(\cwt-\wt_{p,q})$
sharp, i.e.\ $1=\keps(\alpha^v)$ for $v=1,\dots,k$.
\item[(ii)]
It is possible to represent $0\in N_\R$ as a positive linear combination of
the vectors $a^{\alpha^v}\in N$, $v=1,\dots,k$. 
(Recall that $a^\alpha$ is the normal
vector of the supporting hyperplane corresponding to the inequality
``$r_\alpha\geq \keps(\alpha)$''.)
\kenditem
The first property means that, along the chosen walk $\keps$ from 
$p$ to $q$, the arrows $\alpha^1,\dots,\alpha^k$ have all the same direction.
On the other hand, the second property implies that there is a decomposition
$Q_0 = Q_0^+ \sqcup Q_0^-$ with $\alpha^1,\dots,\alpha^k$ pointing from
$Q_0^+$ to $Q_0^-$ and being the only arrows connecting these two parts.
This yields a contradiction.  
\hfill $\Box$
\par


\neu{Uni-4}
Let $\Sigma$ be a complete fan in some $d$-dimensional
vector space $N_\R$ with lattice $N$. 
Denote by $M$ the dual lattice.
By \cite{Ke}, I/\S 3 we know that the cohomology groups of equivariant,
invertible sheaves $\cL$ are $M$-graded and how to calculate their summands
as reduced topological cohomology groups of certain subsets of $N_\R${\DP}
With $r\in M$ and
$A_r:=\{a\in N_\R\, |\; \langle a,r\rangle < \mbox{ord}\,{\cal L}(a)\}$
it follows that
$H^l(X_\Sigma,\cL)_{r}= \widetilde{H}^{l-1}(A_r,k)$ for $l\geq 1$. 
%
%
We would like to use this method to prove a generalization of 
Kodaira-vanishing which holds for toric varieties. 
We restrict the subsets $A_r$ in question to
the $(d-1)$-dimensional unit sphere $S^{d-1} \subset N_{\Rr}$: 
\par
 
{\bf Lemma:}
{\em Let $\phi:N_\R\to\R$ be a continuous function which is linear 
on the cones of $\Sigma$. 
For $B_{r} := \{a \in S^{d-1} \subset N_{\R} \mid 
\langle a,r \rangle < \phi(a) \}$ we denote by
$\Sigma_r\subseteq\Sigma$ the subfan consisting of all cones 
$\sigma\in\Sigma$ such that $\sigma\cap S^{d-1}\subseteq B_r$.
Then, the sets $B_{r}$ and $|\Sigma_r| \cap S^{d-1}$ are
homotopy equivalent.
Moreover, the assertion remains 
true if we replace ``$<$'' by ``$\leq$'' in the definition of $B_{r}$.
}\\
($|\Sigma_r|$ denotes the union of all cones contained in $\Sigma_r$.)
\par

{\bf Proof:}
If the fan $\Sigma$ is simplicial (for instance for smooth varieties
$X_\Sigma$), then it is possible to prove in one strike
that $|\Sigma_r| \cap S^{d-1}$ is a deformation retract of $B_r$. In the
general case, however, it seems to be necessary to project $B_r$ 
successively down dimension by dimension. Using the following fact, each
step can be worked out in the cones of $\Sigma$ separately:\\
{\em
Let $P$ be a compact convex polytope and $H^+$ an open or closed half
space. Then 
$\partial P\cap H^+$ is a deformation retract of $P\cap H^+$.
}\\ 
We leave the proof of this obvious fact to the reader.
\hfill$\Box$
\par


\neu{Uni-5} 
{\bf Proposition:} (Kodaira-vanishing)
{\em
Let $X_\Sigma$ be a complete toric variety with at most Gorenstein
singularities.
Assume that $\cL$ is an equivariant line bundle which is generated by its
global sections.
Then, if the lattice polytope $\Delta\subseteq M_\R$
describing $H^0(X_\Sigma,\,\cL)$ (as explained in \zitat{Mod}{6})
is full-dimensional, we have
$H^l(X_\Sigma,\,\cL\otimes\omega_X)=0$ for $l \geq 1$.
}
\par

{\bf Proof:} 
Denote by $\psi$ the order function of $\cL\otimes\omega_X$.
The order function $\varphi_K$ of the canonical divisor equals $1$ on
the skeleton $\Sigma^{(1)}$; ``Gorenstein'' means that $\varphi_K$ 
is linear on the cones.
Thus, we obtain by the previous Lemma that the sets
$B_r:=\{a\in S^{d-1} \subset N_{\R}\,|\; \langle a,r\rangle < \psi(a)\}$ and
$C_r:=\{a\in S^{d-1} \subset N_{\R}\,|\; \langle a,r\rangle \leq
\psi(a)-\varphi_K(a)\}$ 
are homotopy equivalent.
The first one computes the $r$-th graded piece of
the desired cohomology, and the latter is contractible
since $\psi-\varphi_K$ is the order function of $\cL$, 
i.e.\ $(\psi-\varphi_K)(a)=\mbox{\rm min}\langle a, \Delta\rangle$.
\hfill$\Box$
\par

{\bf Remark:}
The assumption about the dimension of $\Delta$ means that $\cL$ may not be
obtained via pull back from some lower-dimensional variety.
\par


\neu{Uni-6}
{\bf Theorem:}
{\em
Assume $\cwt$ is in general position. Then
$\gExt^{l}_{\cwtM}(\UB,\UB) = 0 \mforall l >0$.
}
\par

{\bf Proof:}
Recall that 
$\gExt^l_{\cwtM}(\UB,\UB)=H^l(\cwtM,\,\underline{\lEnd}\,\UB)=
\bigoplus_{p,q\in Q_0} H^l(\cwtM,\,\UB_{p,q})$.
Then, since $\UB_{p,q}\otimes \omega_{\CM}^{-1} \simeq
\UB_{p,q}\otimes \cL(\cwt)$ is globally generated 
with $\dim\ko{\kD}(\wt_{p,q}+\cwt)=\dim\kD(\cwt)$
(cf.\ Propositions \zitat{Mod}{7} and \zitat{Uni}{3}), the result follows
from Kodaira vanishing. 
\hfill $\Box$
\par


\neu{Uni-7}
The previous theorem asks for the canonical weight to be in general position.
We would like to close this section with a criterion for this fact to
hold. Moreover, we present a criterion for $Q_1 = Q_1(\cwt)$.
\par

{\bf Proposition:}
{\em
The canonical weight $\cwt$ is in general position if and only if
there does not exist any $(t,t)$-wall. 
Moreover, $Q_1(\cwt) = Q_1$ if and only if there are no 
$(1,0)$- or $(1,1)$-walls.}
\par

{\bf Proof: } Assume a $(t^+,t^-)$-wall is given. Then, $\sum_{q \in
Q_0^+}\cwt_q = t^+ - t^-$ by adding the $\kMwt$-equations. Consequently,
$\cwt$ lies on this wall precisely when $t^+ - t^- = 0$. This proves the
first claim.\\
For the second claim, we include the case of a $(1,1)$-wall 
in brackets. Let $Q_0 = Q_0^+ \sqcup Q_0^-$ be the
subdivision defining the wall (see \zitat{Mod}{2}), and denote by $\alpha$ the
unique arrow with $s(\alpha) \in Q_0^+$ and $t(\alpha) \in Q_0^-$ [by
$\beta$ the unique arrow with $s(\beta) \in Q_0^-$ and $t(\beta) \in Q_0^+$].
We show that $\alpha$ is not in $Q_1(\cwt)$ [$\alpha$ and
$\beta$ are not both in $Q_1(\cwt)$].  The set $Q_0^+$ is closed under
successors in $Q \setminus \{ \alpha \}$, but $\sum_{q \in Q_0^+}\cwt_q = 1$
[$ = 0$]. Thus, $Q \setminus \{ \alpha \}$ is not stable [both $Q \setminus
\{\alpha \}$ and $Q \setminus \{ \beta \}$ are not stable]. \\
It remains to show the converse. Assume we have no $(1,0)$-wall and no
$(1,1)$-wall. Thus, $t^+ \geq 2$ for each wall $W$. Assume further that $W$ is a
$(t^+,t^-)$-wall with $t^+ > t^-$. We define the open halfspace $W^+ := \{
\wt \in \Hh \mid \sum_{q \in Q_0^+} \wt_q > 0 \}$, i.e.\
$\cwt \in W^+$. Using the wall crossing
formula from \cite{HilleTor}, Lemma 3.4, we obtain $Q_1(\cwt) = \cup_{\wt \in
\Hh} Q_1(\wt)$. It remains to show that $Q_1 = \cup_{\wt \in
\Hh} Q_1(\wt)$. By assumption, for each $\alpha \in Q_1$ there 
exists a tree in $Q \setminus \{\alpha\}$. 
But for each tree $T$ there exists a weight $\wt$ such
that $T$ is  $\wt$-stable (\cite{HilleThin}, Proposition 2.5). 
This finishes the proof. \hfill $\Box$
\par

%
%

\section{The endomorphism algebra of the universal bundle}\label{End}


\neu{End-1} 
In this section we always assume that $\wt$ is in general position,
i.e.\ the universal bundle on $\wtM$ exists. We
start this section with a result about the endomorphism algebra $\cA$ of 
the universal bundle. This algebra is
non-commutative and finite-dimensional;
in order to formulate the statements in this section, we need some
basic results about those algebras. We denote by $\mrad(\cA)$ the
radical of $\cA$. It consists of 
all strongly nilpotent elements $a$ of $\cA$, that is $(a \cA)^n = 0$ for
$n$ sufficiently large.  Thus $\cA$ is isomorphic to the quotient of
the tensor algebra of the 
$\cA/\mrad(\cA)$-bimodule $\mrad(\cA)/\mrad^2(\cA)$ by some admissible
ideal $I$
$$
\cA \simeq T_{(\cA/\mrad(\cA))}\big(\mrad(\cA)/\mrad^2(\cA)\big)/I.
$$ 
Recall that an ideal $I$ is called {\sl admissible} if
$$
T^n_{(\cA/\mrad(\cA))}(\mrad(\cA)/\mrad^2(\cA)) \subset I \subset
T^2_{(\cA/\mrad(\cA))}(\mrad(\cA)/\mrad^2(\cA))
$$ 
for some $n$. In case $\cA$ is the path algebra of a
quiver, the radical of $\cA$ is the ideal generated by paths of length at
least one. 
\par

A finite-dimensional algebra $\Alg$ is called {\sl basic} if the semisimple
quotient $\cA/\mrad(\cA)$ is a product of fields. Because we deal with basic
algebras over 
an algebraically closed field, this semisimple quotient is a product of copies of
the
ground field. It turns out that each basic finite-dimensional algebra is
isomorphic to the quotient of a path algebra of a finite quiver by some
admissible ideal. Moreover, each finite-dimensional algebra is Morita
equivalent to a basic finite-dimensional algebra, that is the module
categories of both algebras are isomorphic. Thus, if we are interested in
module categories, we may restrict ourselves to modules over basic algebras. \\
The endomorphism ring of $\UB$ is basic precisely when $\UB$
contains only pairwise non-isomorphic direct summands. Consequently, for
$\mEnd(\UB)$ 
to be isomorphic to the path algebra of $Q$, it is necessary that the direct
summands $\UB_q$ are pairwise non-isomorphic. In fact, in the theorem below
we will see that the converse is also true.
\par

 
\neu{End-2}
In any case it would be desirable to
know $\mEnd(\UB)$ 
and its Morita equivalent basic algebra. This leads to the following
definitions (cf.\ \cite{Schofield}):
Let $Q$ be a quiver without oriented cycles. Thus, the path algebra $\ck Q$
is finite-dimensional. 
For an arrow $\alpha\in Q_1$ we define the
{\em localization} $\ck Q[\alpha^{-1}]$ by formally adjoining the
inverse $\alpha^{-1}$ of the arrow $\alpha$, i.e.\ $s(\alpha^{-1}) =
t(\alpha)$, $t(\alpha^{-1}) = s(\alpha)$, and $\alpha^{-1} \alpha =
e_{t(\alpha)}$, $\alpha \alpha^{-1} = e_{s(\alpha)}$ where $e_q$ is the
idempotent in $\ck Q$ corresponding to the vertex $q$. In particular,
$\alpha$ and $\alpha^{-1} \in \ck Q[\alpha^{-1}]$ are {\sl not} 
in the radical. Consequently, $\ck Q[\alpha^{-1}]$ is not basic, because
$\ck Q[\alpha^{-1}]/\mrad(\ck Q[\alpha^{-1}])$ contains the two-by-two full
matrix ring with basis $e_{s(\alpha)},e_{t(\alpha)},\alpha,
\alpha^{-1}$. \\ 
We consider the {\em quotient quiver} $\ko{Q} := Q/(Q_1 \setminus
Q_1(\wt))$ defined by killing the arrows from $Q_1 \setminus Q_1(\wt)$ while
identifying their sources and targets, respectively. Each weight of $Q$
provides in a canonical way a weight of $\ko{Q}$ -- just add the values of
the identified vertices. The corresponding moduli spaces are isomorphic.
The localization $Q[\alpha^{-1}\mid \alpha \notin Q_1(\wt)]$ is Morita
equivalent to $\ck \overline{Q}$. Moreover, this localization  is
finite-dimensional if and 
only if the quotient quiver $\ko{Q}$ contains no oriented cycle. 
\par

{\bf Lemma}:
{\em The line bundle $\UB_q$ is isomorphic to the line
bundle $\UB_p$ if and only if there exists a walk from $p$ to $q$
consisting of arrows {\rm not} in $Q_1(\wt)$.}
\par
 
{\bf Proof: } The class of the bundle $\UB_p^{-1} \otimes \UB_q$ is trivial
if and only if the divisor $\sum_{\alpha \in w}D_\alpha - \sum_{\alpha^{-1}
\in w}D_{\alpha} $ is linearly equivalent to zero for a walk $w$ from $p$
to $q$ in $Q$. Moreover, this divisor is also equivalent to $\sum_{\alpha
\in w \cap Q_1(\wt)}D_\alpha - \sum_{\alpha^{-1} \in w \cap
Q_1(\wt)}D_{\alpha}$, because $D_{\alpha} \sim 0$ for $\alpha \notin
Q_1(\wt)$. This divisor is linearly equivalent to $0$ if and only if $w
\cap Q_1(\wt)$ is a cycle in $Q/(Q_1 \setminus Q_1(\wt))$.
This is true if and only if $w \cap Q_1(\wt) = \emptyset $. \hfill $\Box$
\nopagebreak
\par
\nopagebreak
In particular, the indecomposable direct summands of the universal bundles 
on $\wtM$ and $\cM^{\ko{\wt}}(\ko{Q})$ are
isomorphic; in $\UB(\ko{Q})$ we have just cancelled multiple summands.
\par


\neu{End-3} 
{\bf Theorem:}
{\em The endomorphism algebra $\cA$ of the universal bundle $\UB$ on $\wtM$
is isomorphic to the localization of the path algebra of the quiver by all
arrows {\rm not} in $Q_1(\wt)$. If $Q_1(\wt) = Q_1$, then
$\cA$ is isomorphic to the path algebra of the quiver $Q$.}
\par

{\bf Proof: }
As explained in the previous remarks, we may assume that $Q_1(\wt) = Q_1$.
By Proposition \zitat{Uni}{1} we know that
\[
\mHom_{\wtM}(\UB_p,\UB_q)= H^0(\wtM,\UB_p^{-1}\otimes \UB_q)=
k\cdot\Big\{\mbox{lattice points of } \kD(\wt_{p,q})\Big\}
\]
where $\wt_{p,q}$ is the weight introduced in \zitat{Uni}{2}. 
Since an integral flow in $\kD(\wt_{p,q})$ has values in $\{0,1\}$,
we obtain a bijection between the lattice points of $\kD(\wt_{p,q})$ and
the paths from $p$ to $q$ in $Q$; hence 
$\mEnd_{\wtM}(\UB) \simeq kQ$.
\hfill$\Box$
\par


\neu{End-4} 
Let $\UB$ be a vector bundle on a smooth projective algebraic
variety $\cM$. Let $\cA$ be the endomorphism algebra of $\UB$, which is
finite-dimensional. Moreover, let $\UB = \oplus_{q \in Q_0}\UB_q$ be a
decomposition into  indecomposable direct summands. Then $\cA = \oplus_{q \in
Q_0}e_q\cA$ is a decomposition of $\cA$ into indecomposable projective right
$\cA$-modules. 
We denote by $K^b(\UB_q \mid q
\in Q_0 )$ and $K^b(e_q\cA \mid q \in Q_0 )$ the homotopy category of bounded
complexes $\{ C^i \}$, where each $C^i$ is a direct sum of copies of
$\UB_q$ or $e_q\cA$, respectively. The functor induced by the map $\UB_q
\mapsto e_q\cA$ is an equivalence between $K^b(\UB_q \mid q
\in Q_0 )$ and $K^b(e_q\cA \mid q \in Q_0 )$ because the endomorphism
algebra of $\cA$ viewed as {\sl right} $\cA$-module is $\cA$. 
\par

{\bf Theorem:} {\em Assume $Q$ is a quiver without any $(t,t)$-wall. Then the
equivalence above induces a full and faithful functor 
$$
\cD^b\Big(\mmodstr \cA \Big) \lra \cD^b\Big(\mCoh\big(\cwtM\big)\Big).
$$}
\par

{\bf Proof:} 
We define $p
\leq q$ if $\mHom(\UB_p,\UB_q) \not= 0$. This is a partial order on
$Q_0$ because $\UB_q$ is a line bundle for all $q \in Q_0$.
Consequently, $\mEnd(\UB)$ is a directed algebra (there is an order
on $Q_0$ such that $\mHom_{\cA}(e_p\cA,e_q\cA) = 0$ for $p > q$). A
directed algebra is of finite global dimension, thus the bounded derived
category $\cD^b(\mmodstr \cA)$ of finitely generated right $\cA$-modules is
equivalent to the bounded homotopy category $K^b(e_q\cA \mid q \in Q_0 )$
(cf.\ \cite{Happel} \S 1, 3.3).
Since $\UB$ has no self-extension (Theorem
\zitat{Uni}{6}), the natural functor $K^b(\UB_q \mid q\in Q_0 ) \ra
\cD^b(\mCoh(\cwtM))$ is full and faithful. \hfill $\Box$
\par

{\bf Proof of Theorem \zitat{Int}{3}: }
If there exists no $(1,0)$-wall and no $(1,1)$-wall we have $\cA \simeq \ck
Q$ by Proposition \zitat{Uni}{7} and Theorem \zitat{End}{3}. Thus, the result
follows from Theorem \zitat{End}{4} \hfill $\Box$ 
\par

%
%


{\small
\parbox{9cm}{
Klaus Altmann\\
Institut f\"ur reine Mathematik der\\
Humboldt-Universit\"at zu Berlin\\
Ziegelstr.~13A\\
D-10099 Berlin, Germany\\
E-mail: altmann@mathematik.hu-berlin.de}
\parbox{6cm}{
Lutz Hille\\
Fakult\"at f\"ur Mathematik\\
Technische Universit\"at Chemnitz\\
D-09107 Chemnitz, Germany\\
E-mail: hille@mathematik.tu-chemnitz.de}}

\end{document}